\documentclass{amsart}

\newcommand{\IR}{\mathbb R}
\newcommand{\II}{\mathbb I}
\newcommand{\U}{\mathcal U}
\newcommand{\w}{\omega}
\newcommand{\dens}{\mathrm{dens}}
\newcommand{\lin}{\mathrm{lin}}
\newcommand{\id}{\mathrm{id}}

\newtheorem{theorem}{Theorem}
\newtheorem{problem}{Problem}
\theoremstyle{definition}
\newtheorem{remark}{Remark}

\title[Non-separable sigma-locally compact convex sets]{Topological structure of\\ non-separable sigma-locally compact convex sets}
\author{Iryna Banakh, Taras Banakh, and Katsuhisa Koshino}
\address{I.~Banakh: Institute for Applied Problems of Mechanics and Mathematics of Ukrainian Academy of Sciences, Naukova 3b, Lviv, Ukraine}
\email{ibanakh@yahoo.com}
\address{T.~Banakh: Ivan Franko National University of Lviv (Ukraine) and Jan Kochanowski University in Kielce (Poland)}
\email{t.o.banakh@gmail.com}
\address{K.~Koshino: Doctoral Program in Mathematics, Graduate School of Pure and Applied Sciences, University of Tsukuba, Tsukuba, 305-8571, Japan}
\email{kakoshino@math.tsukuba.ac.jp}

\begin{document}
\begin{abstract} For an infinite cardinal $\kappa$ let $\ell_2(\kappa)$ be the linear hull of the standard othonormal base of the Hilbert space  $\ell_2(\kappa)$ of density $\kappa$. We prove that a non-separable convex subset $X$ of density $\kappa=\dens(X)$ in a locally convex linear metric space is homeomorphic to the space
\begin{itemize}
\item $\ell_2^f(\kappa)$ if and only if $X$ can be written as countable union of finite-dimensional locally compact subspaces;
\item $\II^\w\times \ell_2^f(\kappa)$ if and only if $X$ contains a topological copy of the Hilbert cube $\II^\w$ and $X$ can be written as a countable union of  locally compact subspaces.
\end{itemize}
\end{abstract}
\maketitle

For an infinite cardinal $\kappa$ let $$\ell_2(\kappa)=\big\{(x_i)_{i\in\kappa}\in\IR^\kappa: \sum_{i\in\kappa}|x_i|^2<+\infty\big\}$$ be the standard Hilbert space of density $\kappa$, endowed with the norm $\|(x_i)_{i\in\kappa}\|=\Big(\sum_{i\in\kappa}|x_i|^2\Big)^{1/2}$ and $$\ell_2^f(\kappa)=\big\{(x_i)_{i\in\kappa}\in\ell_2(\kappa):|\{i\in\kappa:x_i\ne 0\}<\infty\big\}$$be the linear hull of the standard orthonormal basis in $\ell_2(\kappa)$. For the smallest infinite cadinal  $\kappa=\w$ the spaces $\ell_2(\w)$ and $\ell_2^f(\w)$ are denoted by $\ell_2$ and $\ell_2^f$, respectively. By $\II=[0,1]$ we shall denote the closed unit interval and by $\II^\w$ the Hilbert cube. A closed subset $A$ of a topological space $X$ is called a {\em $Z$-set} in $X$ if the set $\{f\in C(\II^\w,X):f(\II^\w)\cap A=\emptyset\}$ is dense in the function space $C(\II^\w,X)$ endowed with the compact-open topology. Let $A\subset B$ be two subsets in a linear space $L$. We shall say that $A$ has {\em infinite codimension} in $B$ if the linear hull $\lin(A)$ has infinite codimension in the linear hull $\lin(B)$ of the set $B$.

The following characterization of convex sets homeomorphic to $\ell_2^f$ or $\II^\w\times\ell_2^f$ is obtained due to combined efforts of T.~Dobrowolski \cite{Dob}, D.~Curtis, T.~Dobrowoslki, J.~Mogilski \cite{CDM}, and T.~Banakh \cite{Ban} (see also Theorem 5.3.12 and 5.3.2 in \cite{BRZ}).

\begin{theorem}\label{doban} A convex subset $X$ of a locally convex linear metric space is homeomorphic to the space
\begin{itemize}
\item $l_2^f$ if and only if $X$ is infinite-dimensional and $X$ can be written as a countable union of finite-dimensional compact sets;
\item $\II^\w\times l_2^f$ if $X$ can be written as a countable union of compact $Z$-sets and $X$ contains a subset $Q$ which is homeomorphic to the Hilbert cube and has infinite codimension in $X$.
\end{itemize}
\end{theorem}

In this paper we shall prove a non-separable counterpart of Theorem \ref{doban}. For a topological space $X$ its {\em density} $\dens(X)$ is the smallest cardinality $|D|$ of a dense subset $D\subset X$.

\begin{theorem}\label{main} A non-separable convex subset $X$ of density $\kappa=\dens(X)$ in a locally convex linear metric space is homeomorphic to the space
\begin{itemize}
\item $\ell_2^f(\kappa)$ if and only if $X$ can be written as a countable union of finite-dimensional locally compact spaces;
\item $\II^\w\times\ell_2^f(\kappa)$ if and only if $X$ contains a subspace homeomorphic to the Hilbert cube $\II^\w$ and $X$ can be written as a countable union of locally compact spaces.
\end{itemize}
\end{theorem}

In fact, Theorem~\ref{main} follows from a more general Theorem~\ref{main2} characterizing pairs of convex sets homeomorphic to the pairs $(\ell_2(\kappa),\ell_2^f(\kappa))$ or $(\II^\w\times\ell_2(\kappa),\II^\w\times\ell_2^f(\kappa))$. We say that for topological spaces $A\subset X$ and $B\subset Y$ the pairs $(X,A)$ and $(Y,B)$ {\em are homeomorphic} if there is a homeomorphism $h:X\to Y$ such that $h(A)=B$. By a {\em Fr\'echet space} we understand a locally convex linear complete metric space.

\begin{theorem}\label{main2} Let $X$ be a non-separable convex set of density $\kappa=\dens(X)$ in a Fr\'echet space $L$ and $\bar X$ be the closure of $X$ in $L$. The pair $(\bar X,X)$ is homeomorphic to the pair
\begin{enumerate}
\item $(\ell_2(\kappa),\ell_2^f(\kappa))$ if and only if $X$ can be written as a countable union of finite-dimensional locally compact spaces;
\item $(\II^\w\times\ell_2(\kappa),\II^\w\times\ell_2^f(\kappa))$ if and only if $X$ contains a topological copy of the Hilbert cube $\II^\w$ and $X$ can be written as a countable union of locally compact spaces.
\end{enumerate}
\end{theorem}

Theorem~\ref{main2} will be derived from the following two results, first of which is due to T.~Banakh and R.~Cauty \cite{BC}.

\begin{theorem}[Banakh, Cauty]\label{BC} Each non-separable closed convex set $X$ in a Fr\'echet space  is homeomorphic to the Hilbert space $\ell_2(\kappa)$ of density $\kappa=\dens(X)$.
\end{theorem}

The other result used in the proof of Theorem~\ref{main2} is a topological characterization of the pairs $(\ell_2(\kappa),\ell_2^f(\kappa))$ and $(\II^\w\times\ell_2(\kappa),\II^\w\times\ell_2^f(\kappa))$ due to J.~West \cite{West} (see also \cite{SY} and \cite{Koshino}).

\begin{theorem}[West]\label{t5} A pair $(X,Y)$ of topological spaces $Y\subset X$ is homeomorphic to the pair
$\big(\II^\w\times\ell_2(\kappa),\II^\w\times\ell_2^f(\kappa)\big)$ (resp. $\big(\ell_2(\kappa),\ell_2^f(\kappa)\big)$~) for an infinite cardinal $\kappa$ if and only if
\begin{enumerate}
\item the space $X$ is homeomorphic to $\ell_2(\kappa)$;
\item the space $Y$ can be written as a countable union of (finite-dimensional) locally compact spaces, and
\item the space $Y$ absorbs (finite-dimensional) compact subsets of $X$ in the sense that for each compact (finite-dimensional) subset $K\subset X$, a compact subset $B\subset K\cap Y$, and an open cover $\U$ of $X$ there is a topological embedding $h:K\to Y$ such that $h|B=\mathrm{id}|B$ and $h$ is $\U$-near to the identity embedding $\mathrm{id}:K\to X$.
\end{enumerate}
\end{theorem}

Given a cover $\U$ of a topological space $X$, we say that two maps $f,g:Z\to X$ are {\em $\U$-near} and denote this writing $(f,g)\prec\U$ if for each point $z\in Z$ the doubleton $\{f(z),g(z)\}$ is contained in some set $U\in\U$.
\medskip

\noindent {\it Proof of Theorem~\ref{main2}}. The ``only if'' part in the both statements of Theorem~\ref{main2} are trivial. To prove the ``if'' parts, assume that $X$ is a non-separable convex subset of a Fr\'echet space $L$ and let $\bar X$ be the closure of $X$ in $L$. By Theorem~\ref{BC}, the space $\bar X$ is homeomorphic to the Hilbert space $\ell_2(\kappa)$ of density $\kappa=\dens(\bar X)=\dens(X)$.
Now we consider two cases.
\smallskip

1) Assume that the convex set $X$ can be written as a countable union of finite-dimensional locally compact spaces. By Theorem~\ref{t5}, the homeomorphism of the pairs $(\bar X,X)$ and $(\ell_2(\kappa),\ell_2^f(\kappa))$ will follow as soon as we check that the set $X$ absorbs finite-dimensional compact subsets of $\bar X$. Fix a finite-dimensional compact subset $K\subset \bar X$, a compact subset $B\subset K\cap X$, and an open cover $\U$ of $\bar X$. By the density of $X$ in $\bar X$ and the separability of $K$, there is a separable convex subset $Y\subset X$ of $X$ such that $B\subset Y$ and $K\subset\bar Y$. Moreover, using the fact that $X$ is not separable, we can choose $Y$ so that the closure $\bar Y$ is not locally compact. By Theorem 4.4 of \cite{CDM}, the pair $(\bar Y,\bar Y\cap X)$ is homeomorphic to the pair $(\ell_2,\ell_2^f)$, and by Theorem~\ref{t5}, the set $\bar Y\cap X$ absorbs finite-dimensional compact subsets of $\bar Y$. Consequently, for the finite-dimensional compact subset $K\subset\bar Y\subset\bar X$ there is a topological embedding $f:K\to \bar Y\cap X\subset X$ such that $f|B=\mathrm{id}|B$ and $f$ is $\U$-near to the identity embedding $K\to \bar Y\subset\bar X$. This means that $X$ absorbs finite-dimensional compact subsets of $\bar X$. By Theorem~\ref{t5}, the pair $(\bar X,X)$ is homeomorphic to $\big(\ell_2(\kappa),\ell_2^f(\kappa)\big)$.
\smallskip

2) Next, assume that $X$ contains a subspace $A\subset X$ homeomorphic to the Hilbert cube and $X$ can be written as a countable union of locally compact subsets. By Theorem~\ref{t5}, the homeomorphism of the pairs $(\bar X,X)$ and $(\II^\w\times\ell_2(\kappa),\II^\w\times\ell_2^f(\kappa))$ will follow as soon as we check that the set $X$ absorbs compact subsets of $\bar X$. Fix a compact subset $K\subset \bar X$, a compact subset $B\subset K\cap X$, and an open cover $\U$ of $\bar X$. Using the density of $X$ in $\bar X$ and the separability of the compact set $K\cup A$, we can find a separable convex subset $Y\subset X$ of $X$ such that $A\cup B\subset Y$ and $K\subset\bar Y$. We can assume that $Y=\bar Y\cap X$. Since $X$ is not separable, the compact set $A$ has infinite codimension in $X$. So we can choose $Y$ to be so large that $A$ has infinite codimension in $Y$ and $\bar Y$ is not locally compact. By Proposition~3.1 of \cite{CDM}, each closed locally compact subset of $\bar Y$ is a $Z$-set in $\bar Y$.
It follows that the separable convex set $Y=\bar Y\cap X$ is a countable union of compact $Z$-sets. Since the topological copy $A$ of the Hilbert cube has infinite codimension in $Y$, the convex set $Y$ is homeomorphic to $\II^\w\times \ell_2^f$ by Theorem~\ref{doban}(2). By the Uniqueness Theorem for Skeletoids \cite[Theorem 2.1]{BP}, the pair $(\bar Y,Y)$ is homeomorphic to the pair $(\II^\w\times\ell_2,\II^\w\times\ell_2^f)$ and by Theorem~\ref{t5}, the space $Y$ absorbs compact subsets of $\bar Y$. In particular, for the compact subset $K\subset\bar Y$ there is a topological embedding $f:K\to Y\subset X$ such that $f|B=\id|B$ and $(f,\id)\prec\U$. This means that $X$ absorbs compact subsets of $\bar X$ and the pair $(\bar X,X)$ is homeomorphic to the pair $\big(\II^\w\times\ell_2(\kappa),\II^\w\times\ell_2^f(\kappa)\big)$ according to Theorem~\ref{t5}. This completes the proof of Theorem~\ref{main2}.\hfill\qed
\medskip

\begin{remark} By Proposition~3.2 of \cite{Koshino2}, each infinite-dimensional convex subset $X$ of a Fr\'echet space absorbs finite-dimensional compact subsets of its closure $\bar X$. This fact implies  Theorems~\ref{main}(1) and \ref{main2}(1), see \cite{Koshino2}. It should be mentioned that the paper \cite{Koshino2} contains also a result, a bit weaker than Theorem~\ref{main2}(2). In fact, the current paper appeared as a unification of the results of the first and second authors from one side and the third author from the other side. These results were obtained independently and approximately at the same time, in Spring 2013.
\end{remark}

We do not know if the condition on $Q$ to have infinite codimension in $X$ in Theorem~\ref{doban}(2) can be omitted.

\begin{problem} Assume that a subset $A$ of a Fr\'echet space is homeomorphic to the Hilbert cube $\II^\w$. Does $A$ contain a subset $B$, which is homeomorphic to the Hilbert cube and has infinite codimension in $A$?
\end{problem}

\end{document}